# Asymptotics for a certain group of exponential generating functions

Václav Kotěšovec

e-mail: kotesovec2@gmail.com

Jul 22, 2022

**Abstract:** The exponential generating function for the sequence A143405 in the OEIS is exp(exp(x)*(exp(x) - 1)). This paper analyzes the more general exponential generating function exp(m*exp(b*x) + r*exp(d*x) + s and provides asymptotics for the sequences A143405, A355291, A002872, A002874 and others in the OEIS.

## Main result:

If $m > 0$, $b > d \geq 1$, $r \neq 0$ and $e.g.f. = e^{m\,e^{bx} + r\,e^{dx} + s}$

then $\quad a(n) \sim \dfrac{\left(\dfrac{n}{z}\right)^{n+\frac{1}{2}} e^{m\,e^{bz} + r\,e^{dz} - n + s}}{\sqrt{b\,m\,e^{bz}(bz+1) + d\,r\,e^{dz}(dz+1)}}$

where $\quad z = \dfrac{W\left(\dfrac{n}{m}\right)}{b} - \dfrac{1}{\dfrac{b^2\,m^{d/b}\,n^{1-\frac{d}{b}}\left(W\left(\frac{n}{m}\right)+1\right) W\left(\frac{n}{m}\right)^{\frac{d}{b}-2}}{d\,r} + \dfrac{b}{W\left(\frac{n}{m}\right)} + d}$

and W() is the LambertW function

**In addition**, if $b/d \geq 2$

then $\quad a(n) \sim c * \dfrac{\left(\dfrac{b\,n}{W\left(\frac{n}{m}\right)}\right)^n e^{r\left(\frac{n}{m\,W\left(\frac{n}{m}\right)}\right)^{d/b} + \frac{n}{W\left(\frac{n}{m}\right)} - n + s}}{\sqrt{W\left(\frac{n}{m}\right) + 1}}$

where $\quad c = 1 \quad$ if $\quad b/d > 2$

and $\quad c = e^{-\frac{r^2}{8m}} \quad$ if $\quad b/d = 2$

## Proof:

We will use the Hayman's method, see [1] or [2].

```
In[ ]:= Gz = e^(e^(b x) m + e^(d x) r + s) /. x → z
        az = z * D[Gz, z] / Gz // Simplify
        bz = z * D[az, z] // Simplify
```

Out[ ]=
$$e^{e^{b z} m + e^{d z} r + s}$$

Out[ ]=
$$\left(b\, e^{b z} m + d\, e^{d z} r\right) z$$

Out[ ]=
$$z \left(b\, e^{b z} m + b^2 e^{b z} m\, z + d\, e^{d z} r\, (1 + d\, z)\right)$$

The equation az = 0 cannot be solved, so we determine the asymptotic solution. See also [3]. Since b > d, we determine the main asymptotic term from the equation

```
In[ ]:= Solve[ (b e^(b z) m) z == n, z]
```

Out[ ]=
$$\left\{\left\{z \to \frac{\text{ProductLog}\left[\frac{n}{m}\right]}{b}\right\}\right\}$$

In some cases, the first term of the asymptotic expansion is sufficient, but in most cases at least one additional term must be determined. For this purpose we will use Newton's method, see e.g. [4].

```
In[ ]:= Clear[f];
        f[x_] := (b e^(b x) m + d e^(d x) r) x - n;
        x0 = ProductLog[n/m] / b;
        x1 = x0 - f[x] / D[f[x], x] /. x → x0 // FullSimplify
```

Out[ ]=
$$\frac{n \left(\frac{d^2 e^{\left(-2 + \frac{d}{b}\right) \text{ProductLog}\left[\frac{n}{m}\right]} n\, r}{m^2} + b^2 \left(1 + \text{ProductLog}\left[\frac{n}{m}\right]\right)\right)}{b^3 \left(e^{\text{ProductLog}\left[\frac{n}{m}\right]} m + n\right) + b\, d\, e^{\frac{d\, \text{ProductLog}\left[\frac{n}{m}\right]}{b}} r \left(b + d\, \text{ProductLog}\left[\frac{n}{m}\right]\right)}$$

Now, using the relationship

$$\exp(W(n)) = \frac{n}{W(n)}$$

we get the extended expression for the root of the equation.

*In[ ]:=* **FullSimplify[PowerExpand[**

$$\frac{n\left(\frac{d^2 \left(n/\text{ProductLog}\left[\frac{n}{m}\right]/m\right)^{\left(-2+\frac{d}{b}\right)} n\, r}{m^2} + b^2 \left(1 + \text{ProductLog}\left[\frac{n}{m}\right]\right)\right)}{b^3 \left(n/\text{ProductLog}\left[\frac{n}{m}\right] + n\right) + b\, d\, \left(n/\text{ProductLog}\left[\frac{n}{m}\right]/m\right)^{\frac{d}{b}} r\, \left(b + d\, \text{ProductLog}\left[\frac{n}{m}\right]\right)} - \frac{\text{ProductLog}\left[\frac{n}{m}\right]}{b}\bigg]\bigg]$$

*Out[ ]=*

$$-\frac{1}{d + \frac{b}{\text{ProductLog}\left[\frac{n}{m}\right]} + \frac{b^2\, m^{\frac{d}{b}}\, n^{1-\frac{d}{b}}\, \text{ProductLog}\left[\frac{n}{m}\right]^{-2+\frac{d}{b}}\, \left(1+\text{ProductLog}\left[\frac{n}{m}\right]\right)}{d\, r}}$$

We will now determine the asymptotic formula according to Hayman's method. Since it is an exponential generating function, we must multiply the result by n!, for which we will use the Stirling's formula.

*In[ ]:=* **stirling[n_] := n^n / E^n * Sqrt[2 * Pi * n];**
**Gz / (z^n * Sqrt[2 * Pi * bz]) * stirling[n] // FullSimplify**

*Out[ ]=*

$$\frac{e^{e^{b z}\, m - n + e^{d z}\, r + s}\, n^{\frac{1}{2}+n}\, z^{-n}}{\sqrt{b\, e^{b z}\, m\, z\, (1 + b\, z) + d\, e^{d z}\, r\, z\, (1 + d\, z)}}$$

Now it is time to examine when it is sufficient to use only the first term from the asymptotic solution of the equation and when it is not. For this we will analyze the ratio

*In[ ]:=* **FullSimplify[**

**FullSimplify$\left[e^{e^{b z}\, m + e^{d z}\, r} \Big/ z^n \,/.\, z \to \frac{\text{ProductLog}\left[\frac{n}{m}\right]}{b} - \frac{d\, r\, m^{-\frac{d}{b}}\, \text{ProductLog}\left[\frac{n}{m}\right]^{1-\frac{d}{b}}}{b^2\, n^{1-\frac{d}{b}}}\right] \Big/$**

**FullSimplify$\left[e^{e^{b z}\, m + e^{d z}\, r} \Big/ z^n \,/.\, z \to \frac{\text{ProductLog}\left[\frac{n}{m}\right]}{b}\right]$, {n > 1, b > d, m > 0, d ≥ 1}]**

*Out[ ]=*

$$e^{-e^{\text{ProductLog}\left[\frac{n}{m}\right]}\, m + e^{\text{ProductLog}\left[\frac{n}{m}\right]\left(1-\frac{d\, r\, \left(\frac{m\, \text{ProductLog}\left[\frac{n}{m}\right]}{n}\right)^{-\frac{d}{b}}}{b\, n}\right)}\, m - e^{\frac{d\, \text{ProductLog}\left[\frac{n}{m}\right]}{b}}\, r + e^{\frac{d\, \text{ProductLog}\left[\frac{n}{m}\right]\left(b - \frac{d\, r\, \left(\frac{m\, \text{ProductLog}\left[\frac{n}{m}\right]}{n}\right)^{-\frac{d}{b}}}{n}\right)}{b^2}}\, r}\, \left(1 - \frac{d\, r\, \left(\frac{m\, \text{ProductLog}\left[\frac{n}{m}\right]}{n}\right)^{-\frac{d}{b}}}{b\, n}\right)^{-n}$$

In this complicated expression, we will convert the last term to exponential form and then only deal with the exponent. The expansion is based on series for Log[1-x]

$$\left(1 - \frac{d\,r\left(\frac{m\,\text{ProductLog}\left[\frac{n}{m}\right]}{n}\right)^{-\frac{d}{b}}}{b\,n}\right)^{-n} \sim$$

$$\text{Exp}\left[\frac{d\,r\left(\frac{m\,\text{ProductLog}\left[\frac{n}{m}\right]}{n}\right)^{-\frac{d}{b}}}{b} + \frac{d^2\,r^2\left(\frac{m\,\text{ProductLog}\left[\frac{n}{m}\right]}{n}\right)^{-\frac{2d}{b}}}{2\,b^2\,n} + \frac{d^3\,r^3\left(\frac{m\,\text{ProductLog}\left[\frac{n}{m}\right]}{n}\right)^{-\frac{3d}{b}}}{3\,b^3\,n^2} + \ldots\right]$$

At this point we will further assume that b >= 2d, otherwise the other terms in this expansion will not tends to zero and the asymptotic formula will have to contain one or more next terms, depending on the values of b and d.

*In[ ]:=* `expr = PowerExpand[ExpandAll[-e^{ProductLog[\frac{n}{m}]} m + e^{ProductLog[\frac{n}{m}]\left(1 - \frac{d\,r\left(\frac{m\,\text{ProductLog}[\frac{n}{m}]}{n}\right)^{-\frac{d}{b}}}{b\,n}\right)} m - `

$$e^{\frac{d\,\text{ProductLog}[\frac{n}{m}]}{b}}\,r + e^{\frac{d\,\text{ProductLog}[\frac{n}{m}]\left(b - \frac{d\,r\left(\frac{m\,\text{ProductLog}[\frac{n}{m}]}{n}\right)^{-\frac{d}{b}}}{n}\right)}{b^2}}\,r + \frac{d\,r\left(\frac{m\,\text{ProductLog}[\frac{n}{m}]}{n}\right)^{-\frac{d}{b}}}{b}\Bigg]\Bigg]$$

*Out[ ]=*

$$-e^{\text{ProductLog}[\frac{n}{m}]}\,m + e^{\text{ProductLog}[\frac{n}{m}] - \frac{d\,m^{-\frac{d}{b}}\,n^{-1+\frac{d}{b}}\,r\,\text{ProductLog}[\frac{n}{m}]^{1-\frac{d}{b}}}{b}}\,m -$$

$$e^{\frac{d\,\text{ProductLog}[\frac{n}{m}]}{b}}\,r + e^{\frac{d\,\text{ProductLog}[\frac{n}{m}]}{b} - \frac{d^2\,m^{-\frac{d}{b}}\,n^{-1+\frac{d}{b}}\,r\,\text{ProductLog}[\frac{n}{m}]^{1-\frac{d}{b}}}{b^2}}\,r + \frac{d\,m^{-\frac{d}{b}}\,n^{\frac{d}{b}}\,r\,\text{ProductLog}[\frac{n}{m}]^{-\frac{d}{b}}}{b}$$

While working with limits with Lambert-W function, I discovered a bug in the Mathematica program (since version 12), for more see [5]. I figured out a way to work around this bug. See also [6]. We will use a substitution

$$t = \frac{n}{W(n/m)}$$

$$W(n/m) = \frac{n}{t}$$

$$\exp(W(n/m)) = \frac{n}{m\,W\,(n/m)} = t/m$$

Now the expression transforms to

In[ ]:= $-t/m*m+t/m*e^{-\frac{dm^{-\frac{d}{b}}n^{-1+\frac{d}{b}}r\,\text{ProductLog}\left[\frac{n}{m}\right]^{1-\frac{d}{b}}}{b}}m-(t/m)^{\frac{d}{b}}r+(t/m)^{\frac{d}{b}}*e^{-\frac{d^2m^{-\frac{d}{b}}n^{-1+\frac{d}{b}}r\,\text{ProductLog}\left[\frac{n}{m}\right]^{1-\frac{d}{b}}}{b^2}}r+\frac{dm^{-\frac{d}{b}}n^{\frac{d}{b}}r\,\text{ProductLog}\left[\frac{n}{m}\right]^{-\frac{d}{b}}}{b}\,/.\,\text{ProductLog}\left[\frac{n}{m}\right]\to n/t\,\text{//}\,\text{Simplify}\,\text{//}\,\text{PowerExpand}$

Out[ ]= $\left(-1+e^{-\frac{dm^{-\frac{d}{b}}r\,t^{-1+\frac{d}{b}}}{b}}\right)t+\frac{dm^{-\frac{d}{b}}r\,t^{\frac{d}{b}}}{b}+\left(-1+e^{-\frac{d^2m^{-\frac{d}{b}}r\,t^{-1+\frac{d}{b}}}{b^2}}\right)m^{-\frac{d}{b}}r\,t^{\frac{d}{b}}$

Mathematica already easily handles limits without the Lambert-W function.
If n tends to infinity, t also tends to infinity and we get

In[ ]:= $\text{Assuming}\Big[\{m>0,b>2d,b>0,d>0,\text{Element}[r,\text{Reals}]\},$
$\text{Limit}\Big[\left(-1+e^{-\frac{dm^{-\frac{d}{b}}r\,t^{-1+\frac{d}{b}}}{b}}\right)t+\frac{dm^{-\frac{d}{b}}r\,t^{\frac{d}{b}}}{b}+\left(-1+e^{-\frac{d^2m^{-\frac{d}{b}}r\,t^{-1+\frac{d}{b}}}{b^2}}\right)m^{-\frac{d}{b}}r\,t^{\frac{d}{b}},t\to\text{Infinity}\Big]\Big]$

Out[ ]=
0

In[ ]:= $\text{Assuming}\Big[\{m>0,b>0,d>0,\text{Element}[r,\text{Reals}]\},$
$\text{Limit}\Big[\left(-1+e^{-\frac{dm^{-\frac{d}{b}}r\,t^{-1+\frac{d}{b}}}{b}}\right)t+\frac{dm^{-\frac{d}{b}}r\,t^{\frac{d}{b}}}{b}+\left(-1+e^{-\frac{d^2m^{-\frac{d}{b}}r\,t^{-1+\frac{d}{b}}}{b^2}}\right)m^{-\frac{d}{b}}r\,t^{\frac{d}{b}}\,/.\,b\to 2d,t\to\text{Infinity}\Big]\Big]$

Out[ ]=
$-\frac{r^2}{8m}$

## Conclusion:

In the case b < 2d only one term is not enough and the asymptotics will have a different shape. We have to use the full formula.
For b > 2d, only one term is sufficient. We can use a simplified formula
The most interesting is the case b = 2d, when one term will be sufficient, but the result must be multiplied by a constant factor (we have to consider that we have analyzed the exponent)

In[ ]:= $\exp\left(-\frac{r^2}{8m}\right)$

# Examples:

OEIS A143405, e.g.f.: exp(exp(x)*(exp(x)-1))

*In[ ]:=* $e^{e^{bx}m + e^{dx}r + s}$ /. {m → 1, b → 2, r → -1, d → 1, s → 0}

*Out[ ]=* $e^{-e^x + e^{2x}}$

The asymptotic ratio, 10000 terms. Richardson extrapolation (used a simplified formula), the ratio converges to 1:

*In[ ]:=* `$MaxExtraPrecision = 1000;`

$$\text{funs}[n\_] := A143405[\![n]\!] \Big/ \left( \frac{2^n \, e^{-\frac{1}{8} - n - \sqrt{\frac{n}{\text{ProductLog}[n]}} + \frac{n}{\text{ProductLog}[n]}} \, n^n \, \text{ProductLog}[n]^{-n}}{\sqrt{1 + \text{ProductLog}[n]}} \right);$$

```
Do[Print[N[Sum[(-1)^(m + j) * funs[j * Floor[Length[A143405] / m]] *
    j^(m - 1) / (j - 1)! / (m - j)!, {j, 1, m}], 40]], {m, 100, 1000, 100}]
```

1.009722244529620854463002250172117 66791

1.009342173025063372116874372006202770162

1.009130562842768438104697018837870808640

1.008974293640730129525322587487880136824

1.008859063608003053939673015263185910289

1.008786880229782078206281066776848925354

1.008698671956819945972229510547 67227033

1.008642269271830207799249992405883988959

1.008568754442746352166085145881623635019

1.008512201332924230256235915398874365187

OEIS A355291, e.g.f. exp(exp(x)*(exp(x) + 1) - 2)

In[ ]:= $e^{e^{bx} m + e^{dx} r + s}$ /. {m → 1, b → 2, r → 1, d → 1, s → -2}

Out[ ]=
$e^{-2 + e^x + e^{2x}}$

In[ ]:=
```
$MaxExtraPrecision = 1000;
```

$$\text{funs}[n\_] := A355291[\![n]\!] \bigg/ \left( \frac{2^n\, e^{-\frac{17}{8} - n + \sqrt{\frac{n}{\text{ProductLog}[n]}} + \frac{n}{\text{ProductLog}[n]}}\, n^n\, \text{ProductLog}[n]^{-n}}{\sqrt{1 + \text{ProductLog}[n]}} \right);$$

```
Do[Print[N[Sum[(-1) ^ (m + j) * funs[j * Floor[Length[A355291] / m]] *
    j^(m - 1) / (j - 1) ! / (m - j) !, {j, 1, m}], 40]], {m, 100, 1000, 100}]
```

1.010542379501787494930014959776226709856

1.009935166046433921932389094857063368855

1.009623190471489493313934171345925347703

1.009402644289785880991477112712897534689

1.009244761455312759971557381566459077720

1.009147743554221370276539931830868940432

1.009031041577802719828336573623571052136

1.008957429680662740267900910313196432034

1.008862613679628528022875939034060878209

1.008790503975101954771521362167911665960

OEIS A002872, e.g.f.: exp( (exp(2x) - 3)/2 + exp(x) )

In[ ]:= $e^{e^{bx} m + e^{dx} r + s}$ /. {m → 1 / 2, b → 2, r → 1, d → 1, s → -3 / 2}

Out[ ]=

$e^{-\frac{3}{2} + e^x + \frac{e^{2x}}{2}}$

In[ ]:= `$MaxExtraPrecision = 1000;`

```
funs[n_] := A002872[[n]] / ( (2^n e^{-\frac{7}{4} - n + \sqrt{2} \sqrt{\frac{n}{ProductLog[2n]}} + \frac{n}{ProductLog[2n]}} (\frac{n}{ProductLog[2n]})^n ) / \sqrt{1 + ProductLog[2n]} );
```

```
Do[Print[N[Sum[(-1)^(m + j) * funs[j * Floor[Length[A002872] / m]] *
    j^(m - 1) / (j - 1)! / (m - j)!, {j, 1, m}], 40]], {m, 100, 1000, 100}]
```

1.02000953554461686852152758159674180398 5

1.01892338040793125136785780413714220199 7

1.01836196385222792784664819331514468081 9

1.01796364312220591804940488405372273873 9

1.01767775200663541395381688037301955228 5

1.01750176274810099837454316805535772291 3

1.01728974916468612028259817000084599447 9

1.01715583935809022098184047327646335815 5

1.01698315140826294523511822501069997196 4

1.01685166385989832065745257988297768247 8

OEIS A002874, e.g.f.: exp( (exp(3*x) - 4)/3 + exp(x) )

*In[ ]:=* $e^{e^{bx}m + e^{dx}r + s}$ /. {m → 1/3, b → 3, r → 1, d → 1, s → -4/3}

*Out[ ]=* $e^{-\frac{4}{3} + e^x + \frac{e^{3x}}{3}}$

*In[ ]:=*
```
$MaxExtraPrecision = 1000;

funs[n_] := A002874[[n]] / (3^n e^(-4/3 - n + 3^(1/3) (n/ProductLog[3 n])^(1/3) + n/ProductLog[3 n]) (n/ProductLog[3 n])^n / Sqrt[1 + ProductLog[3 n]]);

Do[Print[N[Sum[(-1)^(m + j) * funs[j * Floor[Length[A002874] / m]] *
    j^(m - 1) / (j - 1)! / (m - j)!, {j, 1, m}], 40]], {m, 100, 1000, 100}]
```

0.9977100498380866662705118292589217275151
0.9981368582100761540845852691045545917726
0.9983452950105909450352310465182708977356
0.9984873306554688327231930924715810821819
0.9985860327312037331349669205654627110148
0.9986453756308877107501731588200956789078
0.9987153706262706570263349497797482720320
0.9987587191853541448727709675499573951234
0.9988136065508655878619321667039649215450
0.9988546191874549299185426142935698076989

See also OEIS A036074, OEIS A036075, OEIS A036076, OEIS A036077, OEIS A036078, OEIS A036079, OEIS A036080, OEIS A036081, OEIS A036082